\documentclass[12pt]{article}

\usepackage{amsmath}
\usepackage{theorem}
\usepackage{amssymb}
\usepackage{latexsym}
\usepackage{a4wide}
\usepackage{graphics}
\newtheorem{thm}{Theorem}
\newtheorem{lem}[thm]{Lemma}
\newtheorem{prop}[thm]{Proposition}
\newtheorem{cor}[thm]{Corollary}

\newtheorem{alg}[thm]{Algorithm}
{\theorembodyfont{\rmfamily} \newtheorem{exa}[thm]{Example}}
\newenvironment{rem}{\noindent{\bf Remark.}}{\newline}
\newenvironment{pf}{\noindent{\bf Proof.}}{\hbox{}\hfill $\Box$}

\newcommand{\Q}{\mathbb{Q}}
\newcommand{\R}{\mathbb{R}}
\newcommand{\C}{\mathbb{C}}
\newcommand{\F}{\mathbb{F}}
\newcommand{\Z}{\mathbb{Z}}

\newcommand{\Gal}{\mathrm{\mathop{Gal}}}
\newcommand{\GL}{\mathrm{\mathop{GL}}}
\newcommand{\End}{\mathrm{\mathop{End}}}
\newcommand{\Lie}{\mathrm{\mathop{Lie}}}
\newcommand{\Tr}{\mathrm{\mathop{Tr}}}

\newcommand{\diag}{\mathrm{\mathop{diag}}}
\newcommand{\lcm}{\mathrm{\mathop{lcm}}}
\newcommand{\rk}{\mathrm{\mathop{rk}}}
\newcommand{\gl}{\mathfrak{\mathop{gl}}}
\newcommand{\g}{\mathfrak{g}}
\newcommand{\gx}{\mathfrak{\mathop{g}}(X)}
\newcommand{\gxf}[1]{\mathfrak{\mathop{g}}_{#1}(X)}

\def\pitem{\advance\leftskip3mm\advance\linewidth-3mm}
\def\mitem{\advance\leftskip-3mm\advance\linewidth3mm}

\begin{document}

\title{Constructing algebraic Lie algebras}
\author{Claus Fieker\\
School of Mathematics and Statistics\\
University of Sydney\\
Australia\\
email: {\tt claus@maths.usyd.edu.au}\and
Willem A. de Graaf\\
Dipartimento di Matematica\\
Universit\`{a} di Trento\\
Italy\\
email: {\tt degraaf@science.unitn.it}}
\date{}
\maketitle

\begin{abstract}
We give an algorithm for constructing the algebraic hull of a given
matrix Lie algebra in characteristic zero.
It is based on an algorithm for finding integral linear dependencies of
the roots of a polynomial, that is probably of independent interest.
\end{abstract}

\section{Introduction}

One of the major tools in the theory of algebraic groups is their
correspondence with Lie algebras. Many problems regarding algebraic groups
can be reformulated in terms of the corresponding Lie algebras, for which they are
generally easier to solve. There is considerable interest in studying
algebraic groups computationally (cf., e.g., \cite{cmt}, 
\cite{grunewald_segal}). Also for this it would be
of great interest to exploit the connection with Lie algebras.
In this paper we treat a question that arises in this context, namely
the problem to decide whether a given Lie algebra corresponds to 
an algebraic group. In particular, a positive solution to this
problems enables us to decide which subalgebras of a Lie algebra
of an algebraic group correspond to algebraic subgroups. To tackle 
this problem we restrict to base fields of characteristic $0$, because
for that case there is a well developed theory of the connection
between algebraic groups and Lie algebras (see \cite{chevii}).
In particular, a connected algebraic group is completely determined
by its Lie algebra.\par 
Let $V$ be a finite-dimensional vector space.
A subgroup $G\subset \GL(V)$ is said to be algebraic if there 
is a set of polynomial functions $P$ on $\End(V)$ such that $G$ consists
of all $g\in \GL(V)$ with $f(g)=0$ for all $f\in P$. To such a group
corresponds a Lie algebra, $\Lie(G)\subset \gl(V)$ (\cite{chevii}, Chapter II, \S 8),
where by $\gl(V)$ we denote the Lie algebra of all endomorphisms of $V$. 
Now a given Lie subalgebra $\g\subset 	\gl(V)$ is called algebraic if there
is an algebraic subgroup $G\subset \GL(V)$ such that $\g = \Lie(G)$.
In \cite{chevii}, Chevalley studied this concept in characteristic $0$,
and gave several sufficient criteria for a $\g\subset \gl(V)$ to be
algebraic. \par
Let $\g\subset \gl(V)$ be any Lie algebra. Then by \cite{chevii},
Chapter II, Theorem 13, there is a unique smallest algebraic Lie
algebra containing $\g$. This algebraic Lie algebra is called the
algebraic hull of $\g$. In this paper we consider the problem of
constructing the algebraic hull for a given $\g\subset \gl(V)$. \par
Based on results of Chevalley we describe an algorithm for constructing the
algebraic hull. The computationally hardest step is to construct the 
splitting field of a polynomial. Since this can be a rather formidable task, we 
subsequently give an algorithm which is similar in nature, but avoids 
the problem of having to construct splitting fields. This is based on
algorithms for finding integral dependencies of algebraic integers. 
Combining complex and $p$-adic approximations to the roots, and the technique
of lattice reduction (LLL), we obtain an algorithm for computing the 
$\Z$-module of integral relations among a given set of algebraic integers. 
In the literature, several somewhat similar methods for solving this
problem are know (cf., e.g., \cite{cohen1} \S2.7.2, \cite{just}). These methods focus
on finding one linear dependency, while our algorithms find (a basis of)
the whole module of linear dependencies. \par
This paper is arranged as follows. In Section \ref{sec:prelim} we introduce
the notation that we use, and summarize a number of results of Chevalley.
Then in Section \ref{sec:basic} we describe the algorithm that makes use
of splitting fields of polynomials. In Section \ref{sec:permmod} we show how
Galois groups can in some instances be of help with constructing the
algebraic hull. This is used in Section \ref{sec:deg4}, where we give
the algebraic hull of the Lie algebra spanned by a semisimple $4\times 4$-matrix.
Then in Section \ref{sec:intdep} algorithms are given for finding integral
linear dependencies among the roots of a polynomial. These algorithms are
then used in Section \ref{sec:algla}, where an algorithm is given for constructing
the algebraic hull of a Lie algebra, avoiding the construction of splitting fields.
Finally, in Section \ref{sec:pract} we report on some practical experiences 
with an implementation of the algorithms in the computer algebra system
{\sc Magma} \cite{magma,magma1}.

\section{Preliminaries}\label{sec:prelim}

Here $F$ will be a field of characteristic $0$. We will use the language of
matrices, rather than that of endomorphisms, as this is more convenient for
calculations. In particular, $\gl(n,F)$ is the Lie algebra of all $n\times n$-matrices
over $F$.
By \cite{chevii}, Chapter II, Theorem 14, a Lie algebra $\g\subset \gl(n,F)$
is algebraic if it is generated by algebraic Lie algebras. It follows that
$\g$ is algebraic if and only if the algebraic hull of the subalgebra spanned
by each basis element of $\g$ is contained in $\g$. Hence
we can compute the algebraic hull of $\g$ if we can compute it in the case
where $\g$ is spanned by one matrix $X$. \par
Let $X\in \gl(n,F)$. Then by $\gxf{F}$ we denote the algebraic hull of the Lie
algebra spanned by $X$. Let $X=S+N$ be the Jordan decomposition of $X$. Then from
\cite{chevii}, Chapter II, Theorem 10 (see also \cite{borel}, \S 7), it follows that
$\gxf{F}=\mathfrak{g}_F(S) \oplus \mathfrak{g}_F(N)$. Moreover, $\mathfrak{g}_F(N)$ is 
spanned by $N$, by \cite{chevii}, Chapter II, \S13, Proposition 1. 
So the problem is reduced to finding $\gxf{F}$ when $X$ is semisimple.\par
The following theorem is proved in \cite{chevii}:
\begin{thm}[Chevalley]\label{thm1}
Let $X\in \gl(n,F)$ be semisimple, and let $K\supset F$ be an algebraic extension
containing the eigenvalues $\alpha_1,\ldots,\alpha_n$ of $X$. Let $U\in \GL(n,K)$
be such that $Y=UXU^{-1}$ is in diagonal form, with the $\alpha_i$ on the diagonal.
Set $\Lambda =\{ (e_1,\ldots,e_n)\in \Z^n\mid \sum_i e_i\alpha_i = 0\}$. 
Then
\begin{enumerate}
\item $\gxf{K} = U^{-1} \mathfrak{g}_{K}(Y)U$ and
$$\mathfrak{g}_{K}(Y)= \{ \diag(a_1,\ldots,a_n) \mid a_i\in K \text{ and }
\sum_i e_ia_i=0 \text{ for all $(e_1,\ldots,e_n)\in \Lambda$} \}.$$
\item $\gxf{F}\otimes K \cong \gxf{K}$.
\item $\gxf{F}\subset A_F(X)$ where $A_F(X)$ is the associative $F$-algebra 
with one generated by $X$. 
\end{enumerate}
\end{thm}

The first part of 1. is straightforward. Let $G_K(X)$ denote the smallest
algebraic subgroup of $\GL(n,K)$ such that its Lie algebra contains $X$.
Then $G_K(X) = U^{-1} G_K(Y) U$ and
$\gxf{K}=\Lie(G_K(X)) = U^{-1} \Lie(G_K(Y)) U = U^{-1} \mathfrak{g}_{K}(Y)U$.
The second part of 1. is \cite{chevii}, \S 13, Proposition 2.
2. follows from the proof of \cite{chevii}, \S 13, Theorem 10. Furthermore,
3. is \cite{chevii}, \S 14, Proposition 14. (There it is shown that $\gxf{F}$ is contained
in the associative algebra (not necessarily with one) generated by $X$. However, for
us it will be more convenient to add the identity.) \par

\section{An algorithm for the algebraic hull}\label{sec:basic}

In this section we use the same notation as in the previous section.
In particular we let $X$ be a semisimple $n\times n$-matrix with coefficients
in the field $F$ of characteristic $0$. We let $K$ be a finite extension of $F$
containing the eigenvalues $\alpha_1,\ldots,\alpha_n$ of $X$. Furthermore,
$\Lambda = \{ (e_1,\ldots,e_n)\in \Z^n \mid \sum_i e_i\alpha_i = 0\}$, and
$\Lambda_\Q = \{ (e_1,\ldots,e_n)\in \Q^n \mid \sum_i e_i\alpha_i = 0\}$.
By $A_F(X)$ we denote the associative algebra with one generated by $X$.
The algorithm is based on the following lemma. 

\begin{lem}\label{lem4}
For $\underline{e}=(e_1,\ldots,e_n)\in \Q^n$ and $i\geq 0$ set $\Delta_i(\underline{e})
= \sum_{k=1}^n e_k \alpha_k^i$. Let $I=X^0,X,\ldots,X^t$ be a basis of $A_F(X)$. Set
$$\Upsilon = \{ (\gamma_0,\ldots,\gamma_t)\in F^{t+1} \mid \sum_{i=0}^t \Delta_i(\underline{e})
\gamma_i = 0 \text{ for all } \underline{e}\in \Lambda_\Q\}.$$
Then $\gxf{F} = \{ \sum_{i=0}^t \gamma_i X^i \mid (\gamma_0,\ldots,\gamma_t)\in \Upsilon\}$. 
\end{lem}

\begin{pf}
Let $Y=\diag(\alpha_1,\ldots,\alpha_n)$. Then there is a $U\in \GL(n,K)$ 
with $UXU^{-1}=Y$. Here $t+1$ is the degree of the minimal polynomial of $X$. Then
since the minimal polynomial of a semisimple matrix is the square free part of 
its characteristic polynomial, the minimal polynomial of $Y$ (over $K$) is the same 
as the minimal polynomial of $X$ (over $F$). Hence $A_K(Y)$ is spanned by $I,Y,Y^2,\ldots,Y^t$.\par
Set $y=\sum_{i=0}^t \gamma_i Y^i$. Write $y(k,k)$ for the entry in $y$ on position $(k,k)$.
Then by Theorem \ref{thm1}, $y\in \mathfrak{g}_{K}(Y)$ if and only if for all 
$\underline{e}\in \Lambda$ we have $\sum_k e_k y(k,k)=0$. It is clear that in this 
statement we may replace $\Lambda$ by $\Lambda_\Q$. Indeed, $\Lambda$ is a subgroup 
of $\Z^n$ and hence it is finitely generated (see, e.g., \cite{schenkman}, Corollary II.3.k). 
Furthermore,  a $\Z$-basis of $\Lambda$ will also be a $\Q$-basis of $\Lambda_\Q$. \par
Now $y(k,k)= \sum_{i=0}^t \gamma_i \alpha_k^i$, and hence $\sum_{k=1}^n e_k y(k,k) = 
\sum_{i=0}^t \Delta_i(\underline{e}) \gamma_i$. Now set 
$$\Upsilon' = \{ (\gamma_0,\ldots,\gamma_t)\in K^{t+1} \mid \sum_{i=0}^t \Delta_i(\underline{e})
\gamma_i = 0 \text{ for all } \underline{e}\in \Lambda_\Q\}.$$
Then by Theorem \ref{thm1} we get that $\mathfrak{g}_K(Y) = 
\{ \sum_{i=0}^t \gamma_i Y^i \mid (\gamma_0,\ldots,\gamma_t)\in \Upsilon'\}$. By the same theorem,
$\gxf{K}= U^{-1} \mathfrak{g}_K(Y) U$ and hence 
$\gxf{K} = \{ \sum_{i=0}^t \gamma_i X^i \mid (\gamma_0,\ldots,\gamma_t)\in \Upsilon'\}$.\par
Now let $X_1,\ldots,X_s$ be any basis of $\gxf{F}$. Then according to 2. of Theorem \ref{thm1} 
they are also a basis of $\gxf{K}$. So $\gxf{K}$ consists of $\sum_i \beta_i X_i$
with $\beta_i\in K$. By 3. of Theorem \ref{thm1}, $X_i\in A_F(X)$.
Hence $\gxf{K}\cap A_F(X)$ consists of $\sum_i \delta_i X_i$ with $\delta_i\in F$. 
We conclude that $\gxf{K} \cap A_F(X) = \gxf{F}$. From this we get the desired conclusion.
\end{pf}

In order to use this result for a practical algorithm we restrict
to the case where $F$ is an algebraic number field. Then the algorithm for 
computing $\gxf{F}$ runs as follows.

\begin{alg}\label{alg0}
Let the notation be as above. We suppose that $F$ is a number field. 
This algorithm computes an $F$-basis for $\gxf{F}$.
\begin{enumerate}
\item Compute an algebraic extension $K\supset F$ containing the eigenvalues
$\alpha_1,\ldots,\alpha_n$ of $X$. 
\item Compute (a $\Q$-basis for) $\Lambda_\Q$.
\item Compute (an $F$-basis for) $\Upsilon$ (where $\Upsilon$ is as 
in Lemma \ref{lem4}). 
\item Return the set consisting of $\sum_{i=0}^t \gamma_i X^i$ where 
  $(\gamma_0,\ldots,\gamma_t)$ runs through the basis of the previous step.
\end{enumerate}
\end{alg}
\begin{pf}
First we show that all steps are computable. First of all, by iteratively factoring 
polynomials over number fields we can compute a number field $K\supset F$ containing 
the eigenvalues $\alpha_1,\ldots,\alpha_n$ of $X$. Furthermore, $K$ has a finite 
$\Q$-basis, and a finite $F$-basis. Then by writing the $\alpha_i$ on a $\Q$-basis of
$K$ we can derive a set of linear equations for $\Lambda_\Q$, hence we can compute a 
basis of this space. Note that $\Delta_i(\underline{e})$ depends linearly on $\underline{e}$.
Hence in order to compute $\Upsilon$ it is enough to consider $\underline e$ in a $\Q$-basis 
of $\Lambda_\Q$. So by writing the $\Delta_i(\underline{e})$ on an $F$-basis of $K$, 
we can derive a set of linear equations for $\Upsilon$. Therefore, we can compute a
basis of this space. The last step is trivially computable.\par
The correctness of the algorithm follows from Lemma \ref{lem4}.
\end{pf}

\begin{exa}\label{exa1}
Let $X\in \gl(4,\Q)$ have minimum polynomial $T^4+bT^2+c$ with $D=b^2-4c$ not 
a square in $\Q$. Then the eigenvalues of $X$ are $\alpha_1=\alpha$,
$\alpha_2=-\alpha$, $\alpha_3=\beta$, $\alpha_4=-\beta$, where 
$\alpha^2 = \frac{1}{2}(-b+\sqrt{D})$ and $\beta^2 = \frac{1}{2}(-b-\sqrt{D})$.
Then $\alpha$ and $\beta$ cannot be proportional over $\Q$ (otherwise
$\alpha^2$ and $\beta^2$ would be as well). Hence the $\alpha_i$ span
a $2$-dimensional subspace of $K$. So $\dim \Lambda =2$, and is spanned
by $\underline{e}^1=(1,1,0,0)$, $\underline{e}^2=(0,0,1,1)$. Then
$\Delta_0(\underline{e}^1)=2$, $\Delta_1(\underline{e}^1)=\Delta_3(\underline{e}^1)=0$,
$\Delta_2(\underline{e}^1)=2\alpha^2$. For $\underline{e}^2$ we get the same
except that $\Delta_2(\underline{e}^2)=2\beta^2$. So
$$\Upsilon = \{ (\gamma_0,\ldots,\gamma_3)\in \Q^3 \mid 2\gamma_0+2\alpha^2\gamma_2 = 
2\gamma_0+2\beta^2\gamma_2=0\}.$$
Hence $\Upsilon$ consists of $(0,\gamma_1,0,\gamma_3)$.
We conclude that $\gx$ is spanned by $X,X^3$.
\end{exa}

\section{The permutation module}\label{sec:permmod}

Here we use the same notation as in the previous section. In this section we make
some observations that on some occasions directly give a basis of $\gxf{F}$.\par 
Let $f$ be the characteristic polynomial of $X$. Let $K$ be the splitting field 
of $f$, and $G=\Gal(K/F)$. We represent $G$ as a permutation
group on the roots $\alpha_1,\ldots,\alpha_n$ of $f$. Let $M$ be the permutation
module of $G$ over $\Q$, i.e., $M$ has basis $w_1,\ldots,w_n$ and $\sigma\cdot w_i
=w_{\sigma(i)}$. On many occasions we will write the elements of $M$ as row vectors.
Then $\sigma(a_1,\ldots,a_n)=(a_{\sigma^{-1}(1)},\ldots,a_{\sigma^{-1}(n)})$.
There is a direct sum decomposition of $G$-modules
$M= M_0 \oplus M_1$, where $M_0=\{ \sum_i a_i w_i \mid \sum_i a_i=0\}$ and 
$M_1$ is spanned by $w_1+\cdots +w_n$. \par
Let $(e_1,\ldots,e_n)\in\Lambda_\Q$ and $\sigma\in G$ then
$0=\sigma( \sum_i e_i\alpha_i) = \sum_i e_i \alpha_{\sigma(i)} = 
\sum_i e_{\sigma^{-1}(i)} \alpha_i$. It follows that $\Lambda_\Q$ is
a $G$-submodule of $M$. So by Maschke's theorem $\Lambda_\Q= V_1\oplus 
\cdots \oplus V_s$, where the $V_r$ are irreducible $G$-submodules.\par
From Lemma \ref{lem4} we recall that $\Delta_i(\underline{e}) = \sum_{k=1}^n
e_k \alpha_k^i$, where $\underline{e}\in \Q^n$. 

\begin{lem}\label{lem1}
Write $f=x^n+a_1 x^{n-1}+\cdots +a_n$. Then the $G$-submodule $M_1\subset M$
occurs in $\Lambda_\Q$ if and only if $a_1=0$. Furthermore, $\Delta_i(\underline{e})
=\Tr(X^i)$, where $\underline{e}=(1,1,\ldots,1)$ spans $M_1$.
\end{lem}

\begin{pf}
We have $a_1=0$ if and only if $\sum_i \alpha_i=0$, hence the first statement.
Set $e=(1,1,\ldots,1)$. Let $Y$ be as in the proof of Lemma \ref{lem4}.
Then $\Delta_i(\underline e)=\sum_k \alpha_k^i =\Tr(Y^i) =\Tr(X^i)$. 
\end{pf}

\begin{lem}\label{lem2}
Suppose that $f$ is square-free and that
$M_0$ is irreducible. Then $a_1=0$ implies $\Lambda_\Q=M_1$
and $a_1\neq 0$ implies $\Lambda_\Q=0$. 
\end{lem}

\begin{pf}
Note that $\Lambda_\Q$ cannot contain $M_0$ since in that case a vector
like $(1,-1,0,\ldots,0)$ would be contained in $\Lambda_\Q$, implying
$\alpha_1=\alpha_2$ (which is impossible because $f$ is square free).
Hence the lemma follows by Lemma \ref{lem1}.
\end{pf}

\begin{cor}\label{cor1}
Suppose that $f$ is irreducible.
Let $A_F(X)$ denote the associative algebra generated by $X$. Suppose
that $G$ is $2$-transitive, or that $F=\Q$ and $n$ is prime. If $\Tr(X)=0$ then
$\gxf{F}$ consists of all $X'\in A_F(X)$ with $\Tr(X')=0$, otherwise $\gxf{F} = A_F(X)$.
\end{cor}

\begin{pf}
If $G$ is $2$-transitive then $M_0$ is irreducible, by 
\cite{james_liebeck}, Corollary 29.10. If $n=p$ is
prime then $M_0$ is irreducible over $\Q$. This can be proved as follows. 
First of all, since $G$ is transitive it contains a $p$-cycle. Now we let
$H$ be the subgroup generated by this $p$-cycle. Then $M$ is also an $H$-module.
Moreover, as $H$-module it is isomorphic to the regular module, i.e., 
to the module afforded by the left action of $H$ on the group algebra
$\Q H$. The $H$-submodules of $\Q H$ are exactly the ideals of $\Q H$.
But $\Q H$ is isomorphic to $\Q[x]/(x^p-1)$, which by the chinese remainder
theorem is isomorphic to $\Q \oplus \Q[x]/(x^{p-1}+x^{p-2}+\cdots +1)$.
We conclude that $\Q H$ splits as the direct sum of two simple ideals.
Hence the $H$-module $M$ is a direct sum of two simple submodules. So the
same holds for the $M$ when viewed as $G$-module. \par
Now the result follows by Lemmas \ref{lem1}, \ref{lem2}.
\end{pf}

In particular, if $G=S_n$ or $G=A_n$ ($n\geq 4)$ then we can easily
compute $\gx$. 

\section{Degree 4}\label{sec:deg4}

Here we use the observations of the previous section to give a complete description
of $\gxf{F}$, where $X$ is a semisimple $4\times 4$-matrix, with irreducible
characteristic polynomial. \par
Let $f=x^4+ax^3+bx^2+cx+d$ be the characteristic polynomial of $X$, and suppose that
it is irreducible. Let $G$ denote the Galois group $\Gal(K/F)$,
where $K$ is the splitting field of $f$. We remark that if $F=\Q$ then it is straightforward
to determine $G$, e.g., by the procedure outlined in \cite{rotman}, Theorem 106.
Note that the case where $G=S_4,A_4$ is settled by Corollary \ref{cor1}.

\begin{prop}\label{prop2}
Suppose that $G$ is not isomorphic to $S_4$ or $A_4$.
Then
\begin{enumerate}
\item if $a=0$ and $a^3-4ab+8c\neq 0$ then $\gxf{F}=\{ X'\in A_F(X)\mid \Tr(X')=0\}$,
\item if $a=0$ and $a^3-4ab+8c = 0$ then $\gxf{F}$ is spanned by $X,X^3$,
\item if $a\neq 0$ and $a^3-4ab+8c\neq 0$, then $\gxf{F}$ is spanned by 
$I,X,X^2,X^3$,
\item if $a\neq 0$ and $a^3-4ab+8c=0$, then $\gxf{F}$ is spanned by 
$I,X,X^2+\frac{4}{3a}X^3$. 
\end{enumerate}  
\end{prop}

\begin{pf}
Since $G$ is a transitive permutation group on $4$ points, not isomorphic 
to $S_4$, $A_4$, there remain the possibilities $G\cong \Z/4\Z$, $G\cong 
D_8$, $G\cong V_4$. These groups have respective generating sets
$\{(1,2,3,4)\}$, $\{(1,2,3,4),(1,3)\}$, $\{(1,2)(3,4), (1,4)(2,3)\}$.
In the first two cases the module $M_0$ decomposes as a direct sum
of two submodules with bases $\{(1,-1,1,-1)\}$, $\{(1,0,-1,0),(0,1,0,-1)\}$
(this holds for both cases). Now $\Lambda_\Q$ cannot contain
the second module (as in that case some roots would be equal).
If $G=V_4$ then $M_0$ decomposes as a direct sum of three submodules,
respectively spanned by $(1,1,-1,-1)$, $(1,-1,1,-1)$, $(1,-1,-1,1)$.
The $G$-module $\Lambda_\Q$ cannot contain two of these vectors,
as otherwise after adding it would follow that two roots are equal.\par
So in all cases, after maybe renumbering the roots, there are the 
following possibilities for $\Lambda_\Q$: $\Lambda_\Q=0$, $\Lambda_\Q$ is
spanned by $(1,1,1,1)$, or by $(1,1,1,1)$, $(1,1,-1,-1)$, or
by $(1,1,-1,-1)$. \par
Let $\alpha_1,\ldots,\alpha_4$ be the roots of $f$. Set 
$a_1=\alpha_1+\alpha_2-\alpha_3-\alpha_4$, $a_2=\alpha_1-\alpha_2
-\alpha_3+\alpha_4$, $a_3=\alpha_1-\alpha_2+\alpha_3-\alpha_4$.
Then the product $a_1a_2a_3$ is a symmetric polynomial in the 
$\alpha_i$, hence can be expressed in terms of the coefficients of
$f$. It turns out that $-a_1a_2a_3=a^3-4ab+8c$. So this number is
zero if and only if $\Lambda$ contains $(1,1,-1,-1)$. This proves 1. 
and 3 (cf. Lemma \ref{lem1}). \par
Suppose that $a^3-4ab+8c=0$. Then we can assume that $\Lambda$ contains
$\underline{e}=(1,1,-1,-1)$. In order to obtain a basis of $\Upsilon$
(cf. Algorithm \ref{alg0}) we have to solve the equation
$\sum_{i=0}^{3} \Delta_i(\underline{e}) \gamma_i=0$. Note that 
$\Delta_0(\underline{e})=\Delta_1(\underline{e})=0$. 
We know that $\alpha_1+\alpha_2-\alpha_3-\alpha_4=0$, and also that 
$\alpha_1+\alpha_2+\alpha_3+\alpha_4=-a$. These two relations are equivalent to 
$\alpha_1+\alpha_2+\frac{1}{2}a=0$ and $\alpha_3+\alpha_4+\frac{1}{2}a=0$. Now 
$\Delta_2(\underline{e})=2\alpha_2^2+a\alpha_2-2\alpha_4^2-a\alpha_4$ as the difference is 
equal to 
$$(\alpha_1-\alpha_2-\frac{1}{2}a)(\alpha_1+\alpha_2+\frac{1}{2}a)+
(-\alpha_3+\alpha_4+\frac{1}{2}a)(\alpha_3+\alpha_4+\frac{1}{2}a) = 0.$$
Similarly, $\Delta_3(\underline{e})=-\frac{3}{4}a(2\alpha_2^2+a\alpha_2-2\alpha_4^2-a\alpha_4)$
as the difference is equal to
\begin{multline*}
(\alpha_1^2-\alpha_1\alpha_2-\frac{1}{2}a\alpha_1+\alpha_2^2+a\alpha_2+\frac{1}{4}a^2)(
\alpha_1+\alpha_2+\frac{1}{2}a)+ \\
(-\alpha_3^2+\alpha_3\alpha_4+\frac{1}{2}a\alpha_3-\alpha_4^2-a\alpha_4-\frac{1}{4}a^2)(
\alpha_3+\alpha_4+\frac{1}{2}a)=0.\end{multline*}
From this it follows that $3a\Delta_2(\underline{e})+4\Delta_3(\underline{e})=0$.
Furthermore, $\Delta_2(\underline{e})=2\alpha_2^2+a\alpha_2-2\alpha_4^2-a\alpha_4=2(\alpha_2-\alpha_4)
(\alpha_2+\alpha_4+\frac{1}{2}a)$. From this we conclude that $\Delta_2(\underline{e})\neq 0$. Indeed,
$\alpha_2-\alpha_4\neq 0$ as $f$ is irreducible. Secondly, $\alpha_2+\alpha_4=-\frac{1}{2}a$
would entail $\alpha_1-\alpha_4=0$ as $\alpha_1+\alpha_2=-\frac{1}{2}a$.\par
Suppose $a\neq 0$. Then the equation $\sum_{i=0}^3 \Delta_i(\underline{e})\gamma_i=0$ is equivalent
to $(-\frac{4}{3a}\gamma_3+\gamma_3)\Delta_3(\underline{e})=0$,  and we have just seen that 
$\Delta_3(\underline{e})\neq 0$.
So $\gamma_3=\frac{4}{3a}\gamma_2$ and 4. is proved. \par
If $a=0$, then $\Delta_3(\underline{e})=0$ and the equation 
$\sum_i \Delta_i(\underline{e})\gamma_i=0$ reduces to
$\gamma_2=0$. Also $\underline{e}'= (1,1,1,1)\in \Lambda$. Then by adding 
$\underline{e}$ and $\underline{e}'$ we see
that $\alpha_2=-\alpha_1$ and $\alpha_4=-\alpha_3$. So $\Delta_1(\underline{e}')=
\Delta_3(\underline{e}')=0$,
and $\Delta_0(\underline{e}')=4$. So we get the equation $4\gamma_0=0$. This proves 2.
\end{pf}

The calculations in the final part of the proof have been done with the help of
{\sc Magma}. Using similar calculations more results of the same flavour can be derived.
Without proof we state the following result.

\begin{prop}
Let $X$ be a $6\times 6$-matrix with irreducible characteristic polynomial 
$p=x^6+ax^5+bx^4+cx^3+dx^2+ex+f$. Suppose that the Galois group has number
4,6,7,8, or 11 in the classification of transitive groups in {\sc Magma}. 
Set 
$$r_1=c+\frac{5}{27}(a^3 - \frac{18}{5}ab) \text{ and } r_2=e- \frac{1}{81}a^5 + 
\frac{1}{27} a^3 b - \frac{1}{3}ad.$$
Then
\begin{enumerate}
\item If $r_1=r_2=0$ then $\gxf{F}$ is spanned by
$$I, X, \frac{a}{2}X^2+X^3, \frac{-5a^3}{54} X^2+\frac{5a}{6} X^4+ X^5.$$
\item If one of $r_1,r_2$ is nonzero then $\gxf{F}$ is equal to $A_F(X)$ if $a\neq 0$,
and equal to $\{X'\in A(X_F) \mid \Tr(X') = 0 \}$ if $a=0$.
\end{enumerate}
\end{prop}

\section{Finding integral dependencies of roots}\label{sec:intdep}

Let $f\in \Q[x]$ be a square free polynomial with roots
$\alpha_1,\ldots,\alpha_n$ in some splitting field $\Gamma \supset \Q$.
Elements of $\Gamma\supset 
       K := \Q[\alpha_1, \ldots, \alpha_n]$ can be represented as
polynomials $g\in \Q[X_1, \ldots, X_n]$ coming from a representation
$K \cong \Q[X_1, \ldots, X_n]/I$ for some zero-dimensional
ideal $I\subset \Q[X_1, \ldots, X_n]$. Although constructive methods for
the construction of $I$ or $K$ are known eg. \cite{yokoyama}, in 
general they are limited to small examples: the splitting field can
have degree as large as $n!$ over $\Q$ and generically, it has.
In what follows we assume $f$ to be monic and integral,
so that $\alpha_i$ are
algebraic integers. We will give algorithms for the following tasks:
\begin{enumerate}
\item Given some $g\in\Z[X_1, \ldots, X_n]$ decide if 
                            $g(\alpha_1, \ldots, \alpha_n)=0$
\item Given $g_j\in \Z[X_1, \ldots, X_n]$ ($1\leq j\leq s$)
  find a $\Z$-module basis for
  $$\Lambda := 
   \{\underline e\in \Z^s \mid \sum_{j=1}^s e_j g_j(\alpha_1, \ldots, \alpha_n) = 0\}.$$
\end{enumerate}
Obviously, both tasks are trivial if exact representations for $K$ or $I$
are known, so we essentially assume that $(K:\Q)$ is too large to allow
direct algebraic constructions to succeed.
Our method will be based on approximate representations of the $\alpha_i$, 
ie. we are going to use the field $\C$ of complex numbers and certain
unramified $p$-adic extensions of $\Q_p$ for our work. For basic 
properties of $p$-adic numbers, we refer to \cite{Serre,Lang}.

Let $p\in \Z$ be a prime number. For any $r\in \Z$, we can write
$r = p^lr'$ for some $r'$ not divisible by $p$. The function
$$v_p : \Z \setminus\{0\}\to \Z : r =p^lr'\mapsto l$$
is called the $p$-adic valuation on $\Z$. We extend $v_p$ to
all of $\Z$ by defining $v_p(0) := \infty$ and extend further to $\Q$
by setting $v_p(a/b) = v_p(a)-v_p(b)$. Via 
 $$|.|_p : \Q \to \Q : x\mapsto p^{-v_p(x)}, 0\mapsto 0$$
 this gives rise to the (normalised) $p$-adic absolute value and 
 thus the $p$-adic topology on $\Q$. The completion $\Q_p$ of $\Q$ wrt 
 to $|.|_p$ is called the field of $p$-adic numbers, it contains the
 $p$-adic integers, the completion $\Z_p$ of $\Z$.

Suppose now that over $\F_p$, the field with $p$-elements, $f$
factors as
$$f = \prod_{i=1}^l f_i \mod p$$
with irreducible, pairwise coprime $f_i\in \F_p[t]$.
Then there is an (unramified) extension $\Gamma/\Q_p$ of degree
$f_p := \lcm_{i=1}^l \deg f_i$ where $f$ splits into linear factors.
Furthermore, there is a unique extension of $|.|_p$ to $\Gamma$ which
is again denoted by $|.|_p$.
Similarly to $\R$ or $\C$, elements in $\Gamma$ cannot, in general,
be represented exactly, instead approximations with a given fixed
precision have to be used.
The advantage of using $\Gamma$ as a splitting field rather than $\C$ or 
$K$ directly, lies in the fact that arithmetic operations in $\Gamma$ incur
less numerical loss of precision that operations with real numbers,
while the algebraic degree of $\Gamma/\Q_p$ is still much smaller than
the degree of $K/\Q$. The main disadvantage of using $\Gamma$ or $\C$ is that,
since there is no exact representation of elements, in general we cannot
decide if an element is zero without additional information.

Lastly, we note that there is exactly one prime ideal $P$ of $\Z_K$ 
(the ring of integers of $K$) such that $\Gamma = K_P$ the 
$P$-adic completion at $P$. For elements $x\in \Z_K$, 
we have $x\in P^k$ if and only if $|x|_p\leq p^{-k}$.

In addition to the $p$-adic information mainly encoded in $\Gamma$ we
are also going to need complex information about elements in $K$.
As a number field $K/\Q$, $K$ admits $(K:\Q)$ many distinct embeddings 
$(.)^{(j)}$ ($1\leq j\leq (K:\Q))$ into the complex numbers.
For any $x\in K$ we define a length:
$$T_2: K\to\R : x \mapsto \sum_{j=1}^{(K:\Q)} |x^{(j)}|^2.$$
Note that $\sqrt{T_2}$ is an Euclidean norm on the $\Q$-vectorspace $K$.
Elementary Galois theory and the inequality between arithmetic and geometric
means can be used to derive non-trivial lower bounds on $T_2(x)$:
\begin{equation}\label{eq2}
\sqrt[(K:\Q)]{N_{K/\Q}(x^2)} \leq \frac 1 {K:\Q} T_2(x)
\end{equation}
which implies for algebraic integers $x\in \Z_K$ that
\begin{equation}\label{eq3}
T_2(x) \geq (K:\Q).
\end{equation}

\begin{rem}
Let $\beta_1$, $\ldots$, $\beta_n\in\C$ be the complex roots of $f$.
In general it is extremely difficult to sort the complex roots in such a way
that $\alpha_i$ corresponds to $\beta_i$ which means that for example
from $\sum_{i=1}^ne_i\alpha_i = 0$ we cannot not, in general, conclude that
$\sum_{i=1}^ne_i\beta_i=0$.\hfill
\end{rem}

After these preliminaries we can now state our algorithm for the
first problem:

\begin{alg}\label{algC}
Let $\alpha_1$, $\ldots$, $\alpha_n\in\Gamma/\Q_p$ 
be the roots of some monic polynomial
$f\in\Z[t]$ and assume that $\Gamma$ is unramified over $\Q_p$.
Set $K := \Q(\alpha_1, \ldots, \alpha_n)$ and
let $g\in \Z[x_1, \ldots, x_n]$ be arbitrary.
This algorithm decides if 
$$g(\alpha_1, \ldots, \alpha_s) = 0.$$
\begin{enumerate}
\item Compute a bound $M>0$ such that 
  $|g(\alpha_1, \ldots, \alpha_n)^{(j)}| \leq M$ for all
  complex embeddings $(.)^{(j)}: K \to \C$.
  Such a bound can be obtained by first computing a bound $M'$
  on the complex roots $\beta_i\in \C$ of $f$ and then estimating
  $|g(\gamma_1, \ldots, \gamma_n)|$ for all choices of
  $|\gamma_i|\leq M'$.
\item Compute a bound $r\geq (K:\Q)$.  
\item Set $k := \lceil \frac r{(\Gamma:\Q_p)}\frac{\log M}{\log p}\rceil$
\item Compute $\tilde\alpha_j$  such that $|\tilde\alpha_j-\alpha_j|_p \leq p^{-k}$
\item Evalute 
  $\tilde G := g(\tilde\alpha_1, \ldots, \tilde\alpha_n)$.
\item If $|\tilde G|_p > p^{-k}$ return {\bf not Zero}, otherwise
      return {\bf IsZero}.
\end{enumerate}
\end{alg}
\begin{pf}
Throughout this proof, we write $\tilde\alpha_i$ for finite precision
approximations to the exact root $\alpha_i\in\Gamma$ that we cannot exactly
represent. Similarly, $G := g(\alpha_1, \ldots, \alpha_n)$
is the exact element that we cannot compute but need to decide if $G=0$ 
and $\tilde G := g(\tilde\alpha_1, \ldots, \tilde\alpha_n)$
is a finite precision approximation.

We first note that since $f\in\Z[t]$ is monic, we have $\alpha_i\in\Z_\Gamma$,
the integral closure of $\Z_p$ in $\Gamma$. Now $g\in\Z[x_1, \ldots, x_n]$
implies $G\in\Z_\Gamma$ as well.
Writing $\tilde\alpha_i = \alpha_i+p^k\beta_i$ with some 
$\beta_i\in \Z_\Gamma$ we obtain from 
the ultrametric property of $|.|_p$:
$$|g(\tilde\alpha_1, \ldots, \tilde\alpha_n)|_p \leq
 |g(\alpha_1, \ldots, \alpha_n)|_p,$$
i.e., there is no precision loss in the evaluation.

Let $K := \Q(\alpha_1, \ldots, \alpha_n)$, as above and
$P$ be the unique prime ideal $\Z_K\supset P|p$ such that
$\Gamma = K_P$.
Then for $x \in \Z_K$ such that $|x|_p\leq p^{-k}$
we obtain $x\in P^k$, thus $N_{K/\Q}(x)\in N_{K/\Q}(P)^k$
and, since $N_{K/\Q}(P)$ is an ideal in $\Z$ generated
by $p^{\Gamma : \Q_p}$:
$$p^{k(\Gamma:\Q_p)}\leq N_{K/\Q}(x).$$

Now, lets assume we have $k$ and $M$ as in the algorithm
and $|\tilde G|_p\leq p^{-k}$.
From
$$\sqrt[(K:\Q)]{N_{K/\Q}(G^2)} \leq \frac {T_2(G)} {K:\Q} \leq  \frac {M^2(K:\Q)} {K:\Q} = M^2$$
we get either $G=\tilde G=0$ or $N_{K/\Q}(G)\geq p^{k(\Gamma:\Q_p)}$.
And thus
$$\frac{k(\Gamma:\Q_p)}{K:\Q}\leq \frac {\log M} {\log p}$$
which contradicts our choices. Thus we conclude, $G=0$ as claimed.
\end{pf}

While the above algorithm can veryify a relation, it does not tell us
how to find one. Also, the precision neccessary to verify relations
can be extremely large, it is essentially linear in $(K:\Q) = \#\Gal(f)$.
In order to use similar ideas to find relations we first need a result
allowing us to get a bound on a basis of the relation lattice:
\begin{thm}\label{masser-bound}
Let $\alpha_1$, $\ldots$, $\alpha_n$ be algebraic integers,
$K := \Q(\alpha_1, \ldots, \alpha_n)$, $r := (K:\Q)$
and define
$$\Lambda := \{\underline e\in\Z^n\mid \sum_{i=1}^ne_i\alpha_i=0\}.$$

Suppose that
$|\alpha_i^{(j)}|\leq M$ for all complex emeddings $(.)^{(j)}: K\to\C$,
then $\Lambda$ has a $\Z$-basis $\underline b_i\in\Z^n$, $1\leq i\leq l$
with 
$$\|\underline b_i\|_\infty\leq n^{n-1} M^{n-1}.$$
\end{thm}
\begin{pf}
The function 
$$f:\Z^n\to \R : \underline e\mapsto \sqrt{T_2(\sum_{i=1}^ne_i\alpha_i)}$$
is a convex distance function in the sense of \cite[p. 250]{Masser}.
Let $m$ be a standart basis element of $\Z^n$, ie.
$m=(m_i)_{1\leq i\leq n}$ and $m_i=0$ for all $i\ne i_0$ while $m_{i_0}=1$.
Then $f(m) = \sqrt{T_2(\alpha_{i_0})}\leq \sqrt r M$.
From \eqref{eq3} we get
for non-zero algebraic integers $x\in K$ that $T_2(x)\geq r$, thus 
$f(m)\geq\sqrt r$, for all $m\in \Z^n$ with $f(m)\neq 0$.
The rest now follows directly from the Proposition of \cite[p. 250]{Masser}.
\end{pf}

The next essential ingredience is the LLL algorithm for lattice reduction.
We need the following property of a reduced basis 
\cite[Theorem 2.6.2.(5)]{cohen1}:
\begin{lem}\label{lll}
Let $\Lambda\subseteq \Z^n$ be a lattice. Suppose that $\Lambda$ contains
linear independent elements $x_1$, $\ldots$, $x_l$ of norm $\|x_i\|_2\leq M$.
Then for a LLL-reduced basis $b_1$, $\ldots$, $b_n$ of $\Lambda$ we have
$\|b_i\|^2_2\leq 2^{n-1}M^2$ for $1\leq i\leq l$.
\end{lem}

Combining the previous results we can now give a first algorithm
for linear dependencies:
\begin{alg}\label{algGen}
Let $f\in\Z[t]$ be monic and $\alpha_1$, $\ldots$, $\alpha_n\in\Gamma/\Q_p$
be the roots of $f$ in some unramified extension of $\Q_p$ of degree $f_p$.
We assume that elements in $\Gamma$ are represented as vectors in
$\Q_p^{f_p}$ wrt. some fixed basis $\omega_1$, $\ldots$, $\omega_{f_p}$.
Furthermore, let $g_i\in \Z[x_1, \ldots, x_n]$ be arbitrary ($1\leq i\leq s$)
and define
$$\Lambda := \{\underline e\in \Z^s \mid \sum_{i=1}^s e_ig_i(\alpha_1, \ldots, \alpha_n)=0\}.$$
This algorithm computes a $\Z$-basis for $\Lambda$.
\begin{enumerate}
\item Compute a bound $M'>0$ such that for each $i$ we have
  $|g_i(\alpha_1, \ldots, \alpha_n)^{(j)}|\leq M'$, for all
  complex embeddings $(.)^{(j)}: K \to \C$.
\item Set $N := s^{s-1}M^{s-1}$
\item Set 
  $k := \lceil \frac{(\Q(\alpha_1, \ldots, \alpha_n):\Q)}{f_p}\frac{\log NMs}{\log p}\rceil$
\item Set $\lambda := N^2 2^{s-1}$
\item Compute $\tilde\alpha_j$  such that $|\tilde\alpha_j-\alpha_j|_p \leq p^{-k}$
\item Compute $\tilde\beta_i := g(\tilde\alpha_1, \ldots, \tilde\alpha_n)$
  for $1\leq i\leq s$,
  and form a matrix $B$ where the $i$th row contains the lift of coefficients
  of $\tilde\beta_i$ as elements to $\Z$.
\item Form a big matrix $\tilde B\in\Z^{(s+f_p)\times (s+f_p)}$ 
  by first concatenating $I_s$ and $\lambda B$
  to get $(I_s | \lambda B)$ and then appending a matrix
  $(0I_s | \lambda p^kI_{f_p})$ to the bottom.
\item Apply the LLL algorithm to the rows of $\tilde B$ obtaining a new matrix
  $L = (L_{i,j})_{1\leq i,j\leq f_p+s}$.
\item The lattice $\Lambda$ is generated by 
  $(L_{i,j})_{1\leq i\leq l, 1\leq j\leq s}$ where $l$ is the index
  of the last row $L_i$ of $L$ with norm $\|L_i\|_2<\lambda$.
\end{enumerate}
\end{alg}
\begin{pf}
Using Theorem \ref{masser-bound} we see that $N$ is a bound for the 
maximum norm of a length of a basis-relation, so that $NMs$ is a
bound for the complex embedding $|(.)^{(j)}|$ of a possible relation.
The precision is now
chosen in the same way as in Algorithm \ref{algC} so that a possible relation
$\underline e\in \Z^s$ with $\|e\|_2<N$ and 
$|\sum_{i=1}^se_ig_i(\alpha_1, \ldots, \alpha_n)|_p<p^{-k}$ has to be zero.

In the matrix $L$, the $s$-leftmost columns encode the transformations
applied to $B$ while the rightmost columns give the evaluated relation:
$$\lambda\sum_{i=1}^s L_{j,i}g_i(\tilde\alpha_1, \ldots, \tilde\alpha_n) 
  = \sum_{i=1}^{f_p} L_{j,i+s}\tilde\omega_i + p^kx$$
(for some $x\in \Z_\Gamma$).  So we see that if there is a relation
$\sum_{i=1}^s e_i g_i(\alpha_1,\ldots,\alpha_n)=0$, then the $\Z$-span of the
first $s$ rows of $\tilde{B}$ contains a vector 
$(e_1,\ldots,e_n,u_1,\ldots,u_{f_p})$, with $u_i\in \lambda p^k\Z$. So by adding
suitable multiples of the last $f_p$ rows of $\tilde{B}$ we get that 
$(e_1,\ldots,e_n,0,\ldots,0)$ lies in the $\Z$-span of the rows of $\tilde{B}$.

Our choice of $k$ and $\lambda$ now ensures the following facts:
\begin{enumerate}
\item If and only if $(L_{j, i})_{1\leq j\leq s}$ is a true relation,
      $(L_{j, i+s})_{1\leq j\leq f}$ is zero. 
\item If $(L_{j, i+s})_{1\leq j\leq f}$ is not zero, then
  $\|(L_{i,j})_{1\leq j\leq s+f}\|_2\geq \lambda >N$
\item If there are relations within the bounds of Theorem \ref{masser-bound},
  then the LLL will find them since by Lemma \ref{lll}, there
  must be rows in $L$ with norm bounded by $2^{s-1}N <\lambda$
  which implies that they are relations.
\end{enumerate}
\end{pf}

In applying the above algorithm, the main problem is the huge precision
$k$ needed to guarantee correctness. Since the precision determines
directly the bit-length of the entries of $B$, it is the crucial parameter
for the runtime of the LLL algorithm. By \cite{stehle} we know that
the runtime depends quadratically on the bit-length on the input, thus
we need to try to reduce the precision.
Since verification of a relation (using Algorithm \ref{algC}) is computationally
much easier than finding a relation, one method is to just use the above
Algorithm \ref{algGen} with a smaller precision, say
$\lceil 1.5(\log N/\log p)\rceil$, apply the LLL alorithm and test the 
relations obtained. In case Algorithm \ref{algC} fails to verify a relation
obtained this way, we increase the precision and try again. The proof 
of the correctness shows that this method must terminate with the 
correct answer.

If the Galois-action on the $p$-adic roots $\alpha_1$, $\ldots$, $\alpha_n$
is known, then we can improve the runtime substantially, by using
the fact that if $\sum_{i=1}^sg_i(\alpha_1, \ldots, \alpha_n) = 0$ then
we also have
$\sum_{i=1}^sg_i(\alpha_{\sigma 1}, \ldots, \alpha_{\sigma n}) = 0$
for all $\sigma\in G=\Gal(f)$. This allows us to replace LLL by much
faster echelon algorithms over $\Z/p^k\Z$ followed by rational
reconstruction.
\begin{alg}\label{algGal}
Let $f\in\Z[t]$ be monic and $\alpha_1$, $\ldots$, $\alpha_n\in\Gamma/\Q_p$
be the roots of $f$ in some unramified extension of $\Q_p$ of degree $p$.
Furthermore, let $G = \Gal(f) < S_n$ be given explicitly, ie
$\sigma \alpha_i = \alpha_{\sigma i}$.
Now, let $g_i\in \Z[x_1, \ldots, x_n]$ be arbitrary ($1\leq i\leq s\leq \#G$)
and define
$$\Lambda := \{\underline e\in \Z^s \mid \sum_{i=1}^s e_ig_i(\alpha_1, \ldots, \alpha_n)=0\}.$$
This algorithm computes a $\Z$-basis for $\Lambda$.
\begin{enumerate}
\item Compute a bound $M>0$ such that for each $i$ we have
  $|g_i(\beta_1, \ldots, \beta_n)^{(j)}|<M'$ 
\item Set $N := s^{s-1}M^{s-1}$
\item Set 
  $k := \lceil 2\log NM/\log p\rceil$
\item \label{alg:s1} 
  Select a set $S\subseteq G$ of size $s$, containing
  the identity of $G$.
\item repeat\label{repeat}
\item \pitem Compute $\tilde\alpha_j$  such that $|\tilde\alpha_j-\alpha_j|_p \leq p^{-k}$
  \item Set $B :=()$ a matrix with $s$ rows and $0$ columns. 
  \item \label{alg:s2}for $\sigma\in S$ do

  \item \pitem Compute $\tilde\beta_i := 
                g_i(\tilde\alpha_{\sigma 1} \ldots, \tilde\alpha_{\sigma n})$
      for $1\leq i\leq s$,
      and form a matrix $\tilde B$ where the $i$th row 
      contains the lift of coefficients
      of $\tilde\beta_i$ as elements to $\Z$.
    \item \label{alg:s3}
      Set $B := (B|\lambda\tilde B)$, ie. append $\lambda\tilde B$ to the right of $B$.

  \item \mitem\label{stepN}
    Apply HNF techniques to compute a the nullspace $N$
    of $B\in (\Z/p^k\Z)^{s\times s}$ in echelon form. 
  \item \label{alg:s4}Use rational reconstruction to find (if possible) the unique
   $\tilde N\in \Q^{l\times s}$ such that $\tilde N \cong N \mod p^k$.  If
   this fails, increase the set $S$ by randomly selecting at most $0.2\#S$
   elements in $G\setminus S$ and $k := \lceil 1.2 k\rceil$ and go back to
   step \ref{repeat}.

\item \label{alg:s6}Compute a matrix $S\in \Z^{l\times n}$ such that
  $S$ is a $\Z$-basis for the intersection of the $\Q$-vectorspace
  with basis $\tilde N$ and $\Z^s$ (using some saturation method)
\item \label{alg:s7}Apply the LLL algorithm to $S$ to obtain a LLL reduced basis $L$.
\item \label{alg:s5}Set $k := \lceil 1.2k\rceil$ and increase the set $S$ by randomly
selecting at most $0.2\#S$ elements in $G\setminus S$.  
\item \mitem until\label{stepVerify}
all rows $L_i$ of $L$ are norm bounded: $\|L_i\|_2<N$ 
  and are true relations by Algorithm \ref{algC}
\end{enumerate}
\end{alg}
\begin{pf}
Let $K := \Q(\alpha_1, \ldots, \alpha_n)$, then, since $\Gamma$ is a splitting
field for $f$, we have
$K \otimes_\Q \Gamma \cong \Gamma^G = \Gamma^{(K:\Q)}$ and the embedding is
given via
$\alpha_i \mapsto (\sigma(\alpha_i))_{\sigma\in G}$.
Furthermore, as $\Q_p$-vectorspace we have $\Gamma \cong \Q_p^f$ and
$K\otimes_\Q \Gamma \cong \Q_p^{f(K:\Q)}$ and the embedding extends via
composition with $\Gamma\ni\gamma = \sum_{i=1}^f\gamma_i\omega_i \mapsto (\gamma_i)_{1\leq i\leq f}\in \Q_p^f$.
If we apply this embedding to
$V := [g_1(\underline\alpha), \ldots, g_s(\underline\alpha)]_\Q$ we
get $V\otimes\Gamma = [(g_1(\sigma\underline\alpha))_{\sigma\in G}, \ldots
 (g_s(\sigma\underline\alpha))_{\sigma\in G}]_\Gamma$.
The fact that extensions of scalars preseves dimension, 
$$\dim_\Q(V) = \dim_\Gamma(V\otimes \Gamma) = \dim_{\Q_p}(V\otimes \Gamma)$$
implies that eventually $\rk_{\Q_p} B = \dim_\Q(V)$.
Similarly, the increasing precision ensures that the rational reconstruction
will be successful, eventually. To be more precise:
Assume that $S$ is large enough so that $\dim_\Q(V) = \rk_{\Q_p} B$
and let $M\in\Q^{l\times s}$ be a $\Q$-basis matrix for the $\Q$-nullspace
of $B$. Without loss of generality, we can furthermore 
assume that $M=(M_{i,j})$ is in echelon form, thus $M$ is uniquely 
defined by $V$.
If the precision $k$ is chosen so that $p^k > \max h(M_{i,j})^2$
for the naive height $h:\Q \to\Z : p/q \mapsto \max (p, q)$, then
$N$ in step \ref{stepN} will allow to compute $M$ by reconstruction.

If however $S$ is too small, it may happen that $\rk_{\Q_p} B<\dim_\Q(V)$
which means that the matrix $N$ computed in step \ref{stepN} does not
represent an approximation to $M$. In this case, either the reconstruction
fails or the reconstructed relations cannot be verified in step 
\ref{stepVerify}.
\end{pf}
%%
% it is neccessary to increase the size of S:
% (Jasper): \sqrt a, \sqrt b
%           -\sqrt a, -\sqrt b
% for K = Q(\sqrt a, sqrt b) and S = {id, ...}
% has rk 1, but Q-dim(V) = 2.
% I think any set S of size 2 will fail...

In steps \ref{alg:s1}, \ref{alg:s4} and \ref{alg:s5} we increase the size of
the matrix by exploiting the Galois action. This is done hoping that
generically the new columns are independent from the previous ones as to
increase the rank of $B$. It is clear that if $S=G$ then the 
$\Z$-nullspace of $B$
would be precisely $\Lambda$, however if $S=\{I_G\}$ then there will always be
spurious relations, some of which may be small in size.

Also, using well known techniques that generalize rational reconstruction
to number fields \cite{fieker,belabas,geissler} and omitting steps 
\ref{alg:s6} and \ref{alg:s7} the algorithms can easily be extended
to find $R$-relations instead of $\Z$-relations for $R$ being any order
in some number field.

\begin{rem}
If the Galois action is known then on some occasions we can use 
Algorithm \ref{algC} to give a more efficient 
algorithm for finding a basis of $\Lambda = \{(e_1,\ldots,e_n)\mid \sum_i 
e_i\alpha_i = 0\}$. As in Section \ref{sec:permmod} we denote 
the permutation module of $G$ by $M$ (where $G$ is the Galois group). 
We assume that $M$ has a unique decomposition as direct sum of
irreducible $G$-modules, $M=V_1\oplus \cdots \oplus V_r$. The uniqueness of
this decomposition is equivalent to all of the $V_i$ being non-isomorphic.
In that case we can compute the direct sum decomposition by computing a 
maximal set of orthogonal primitive idempotents in the centre of the
algebra 
$$\mathrm{End}_G(M) = \{ T: M\to M \mid T \text{ is linear and } T(\sigma(v)) = \sigma(T(v)) 
\text{ for $v\in M$ and $\sigma\in G$} \}.$$
It also follows that $\Lambda_\Q = V_{i_1} \oplus \cdots\oplus V_{i_k}$. 
Now for each $V_i$ we do the following. For each element $(e_1,\ldots,e_n)$ in
a basis of $V_i$ we check whether $\sum_i e_i\alpha_i=0$ with Algorithm
\ref{algC}. Then $\Lambda$ is
equal to the direct sum of the $V_i$ that pass this test.\hfill
\end{rem}

\section{A second algorithm for the algebraic hull}\label{sec:algla}

In this section we assume that the base field $F$ is $\Q$. We describe a second 
algorithm for constructing $\gxf{\Q}$, for a semisimple
matrix $X$. It is similar in nature to the algorithm of Section \ref{sec:basic}.
But instead of constructing the splitting field we make use of the algorithms
in the previous section. For simplicity we assume that the characteristic
polynomial is square free. The generalisation to the general case is
straightforward. We use the same notation as in Section \ref{sec:basic}.\par
First we find $\Lambda := \{\underline e\in\Z^n | \sum_{i=1}^n e_i\alpha_i=0\}$
using $g_i := x_i$ and either
Algorithm \ref{algGen} or Algorithm \ref{algGal}.
Let $\underline{e}=(e_1,\ldots,e_n)$ be a basis element
of $\Lambda$. The second step consists in solving the equations that define $\Upsilon$
(cf., Lemma \ref{lem4}). For $i\geq 0$ set $g_i(\underline{e}) =
\sum_{k=1}^n e_k x_k^i$, and $\Delta_i(\underline{e}) := 
g_i(e)(\alpha_1, \ldots, \alpha_n)$. Let $t+1$ be the degree of the minimal polynomial
of $X$. Then again with Algorithm \ref{algGen} or \ref{algGal}
we find all integral (or equivalently, rational) linear
dependencies of the $\Delta_i(\underline{e})$, i.e. all vectors $u=(u_1,\ldots,u_t)\in \Q^t$
with $\sum_i u_i \Delta_i(\underline{e}) = 0$. Let $M(\underline{e})$ 
denote the $\Q$-vector space spanned by all those vectors $u$. Then
$\Upsilon$ is equal to the intersection of all $M(\underline{e})$, where 
$\underline{e}$ runs through a basis of $\Lambda$. So this way we find a basis
of $\Upsilon$, and hence a basis of $\gxf{\Q}$ (cf. Algorithm \ref{alg0}).

\section{Examples}\label{sec:pract}
To generate a set of input examples we used the database of 
polynomials over the rationals with given Galois groups by
Kl\"uners and Malle \cite{database}. In this database the $n$-th
transitive permutation group on $d$ points is denoted ${}_dT_n$.
For each polynomial of degree
$d$ ($6\leq d\leq 12$) with Galois group isomorphic to ${}_dT_n$ we
computed the companion matrix $X$ of $f$ and used this as input to 
our algorithms.
In Figure \ref{figP} we plot the running times for the computation
of $\gxf{\Q}$ using both the algorithm in Section \ref{sec:basic} with an 
exact, algebraic representation of the splitting field of $f$ as well
as the algorithms in Section \ref{sec:algla} against the logarithm of
the group size. From the data presented, it is clear that the
runtime of all three algorithms depends mainly on the size of the
Galois group of $f$, ie. the degree of the splitting field. Also,
clearly, the algebraic representation of the splitting field has the worst
runtime behaviour.
In the second figure, we use a variation of the algorithms in Section
\ref{sec:algla}, where instead of using the bounds from
Algorithm \ref{algC}, we compute the relations with a much smaller bound
and ``verify'' them using twice the $p$-adic precision. While this
of course does not give guaranteed results, nevertheless, in all
cases where the bounds were small enough to use them, the output obtained
thus was correct. Since this approach does not directly depend
on the size of the splitting field, we can use this for larger degrees.

{}From both the figures we notice that for the purpose of computing
algebraic hulls, it does not matter if Algorithm \ref{algGen} or \ref{algGal}
is used. For proven results, the time is always dominated by the proof step
while the actual computation takes only negligable time - even in large degrees
and large Galois groups.

\begin{figure}
%\epsfbox{allP.0}
\includegraphics{figure1.ps}
\caption{Time vs. $\log \#\Gal(f)$ for $f$ of degree $6, 8, 9$ and $10$
and all transitive groups with proven bounds. {\bf C} is used for data
coming from the algebraic, exact representation of the splitting field, 
{\bf B} is time using \ref{algGal} and {\bf A} is using \ref{algGen}.}
\label{figP}
\end{figure}

\begin{figure}
%\epsfbox{allNp.0}
\includegraphics{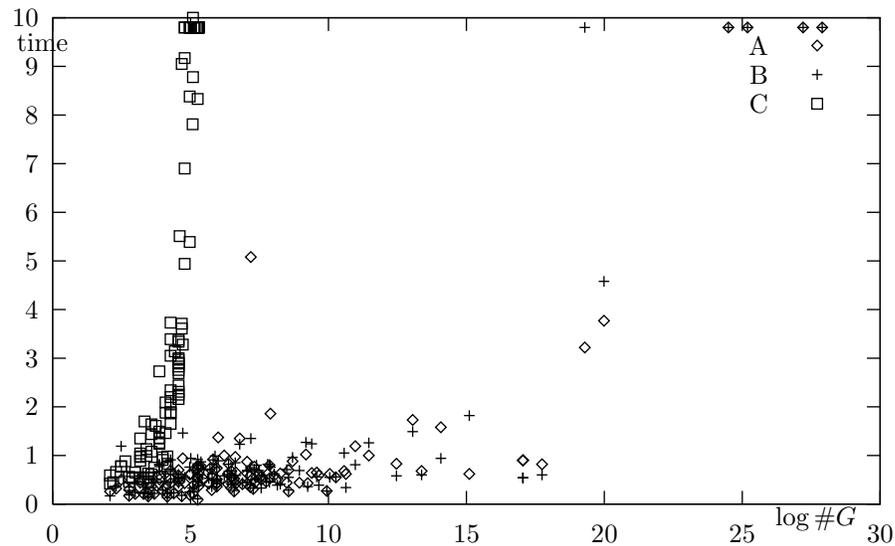}
\caption{Time vs. $\log \#\Gal(f)$ for $f$ of degree $6, 8, 9, 10, 12, 14$ and $15$
and all transitive groups, using heuristic bounds. {\bf C} is used for data
coming from the algebraic, exact representation of the splitting field, 
{\bf B} is time using \ref{algGal} and {\bf A} is using \ref{algGen}.}
\label{figNp}
\end{figure}

\end{document}